# Introduction

Holomorphic dynamics posed several closely related key problems going back to Fatou and Julia: density of axiom A maps, local connectivity of the Julia and Mandelbrot sets, measure and dimension of the above sets. A great deal of progress has been achieved in these problems during the last decade, but they are still in the focus of modern reasearch. Most of the discussion in the chapter is concentrated on these problems. We will not quote any particular results: the reader will find references and a variety of viewpoints on the subject inside the articles.

For a general introduction to holomorphic dynamics we recommend one of the following surveys: [**Be**], [**Bl**], [**C**], [**EL**], [**L**] and [**M**].

The chapter is organized into five sections which cover a good part of the field:

1 Quasiconformal Surgery and Deformations
2 Geometry of Julia Sets
3 Measurable Dynamics
4 Iterates of Entire Functions
5 Newton's Method

The chapter partially arose from an earlier problem list [**Bi**]. This preprint will be published by Springer-Verlag as a chapter in *Linear and Complex Analysis Problem Book* (eds. V. P. Havin and N. K. Nikolskii). We hope it will help to fix up the present state of affairs of the field and to stimulate further development. We thank everybody who has contributed to the chapter. We would also like to thank M. Herman and J. Milnor for many helpful comments.

<div align="center">Ben Bielefeld     Mikhail Lyubich</div>

# Table of Contents

It is possible to investigate rational functions using the technique of quasiconformal surgery as developed in [**DH**2], [**BD**] and [**S**]. There are various methods of gluing together polynomials via quasiconformal surgery to make new polynomials or rational functions. The idea of quasiconformal surgery is to cut and paste the dynamical spaces for two polynomials so as to end up with a branched map whose dynamics combines the dynamics of the two polynomials. One then tries to find a conformal structure that is preserved under this branched map of the sphere to itself, so that using the Ahlfors-Bers theorem the map is conjugate to a rational function. There are several topological surgeries which experimentally seem to exist, but for which no one has yet been able to find a preserved complex structure.

The first such kind of topological surgery is **mating** of two monic polynomials with the same degree. (Compare [**TL**].) The first step is to think of each polynomial as a map on a closed disk by thinking of infinity as a circle worth of points, one point for each angular direction. The obvious extension of the polynomial at the circle at infinity is $\theta \mapsto d\theta$ where $d$ is the degree of the polynomial. Now glue two such polynomials together at the circles at infinity by mapping the $\theta$ of the first polynomial to $-\theta$ in the second. Finally, we must shrink each of the external rays for the two polynomials to a single point. The result should be conjugate to a rational map of degree $d$. (Surprisingly this construction sometimes seems to make sense even when the filled Julia sets for both polynomials have vacuous interior.)

For instance we can take the rabbit to be the first polynomial, that is $z^2 + c$ where the critical point is periodic of period 3 ($c \sim -.122561 + .744862i$). The Julia set appears in the following picture.

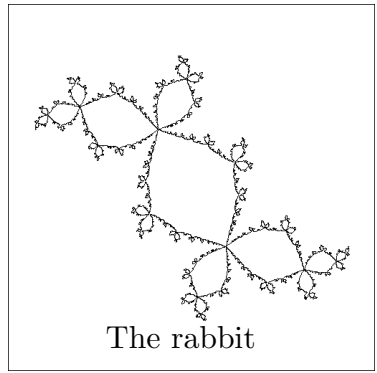

The rabbit

Then for the second polynomial we could take the basilica, that is $z^2 - 1$ (it is named after the Basilica San Marco in Venice. One can see the basilica on top and its reflection in the water below). The Julia set for the basilica appears in the following figure.

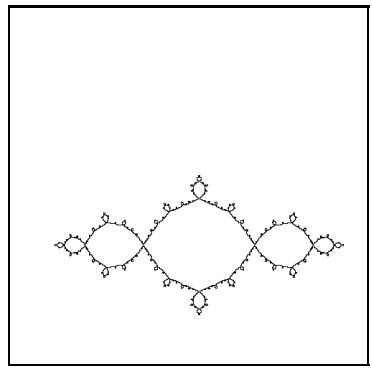

The basilica

Next we show the basilica inside-out $(\frac{z^2}{z^2-1})$ which is what we will glue to the rabbit.

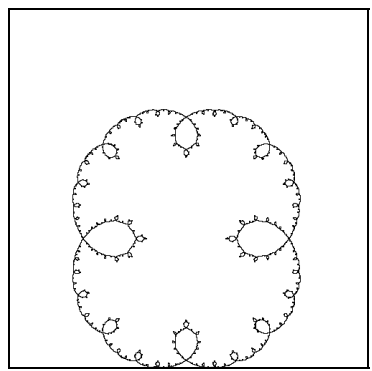

The inside-out basilica

And finally we have the Julia set for the mating $(\frac{z^2+c}{z^2-1}$ where $c = \frac{1+\sqrt{-3}}{2})$.

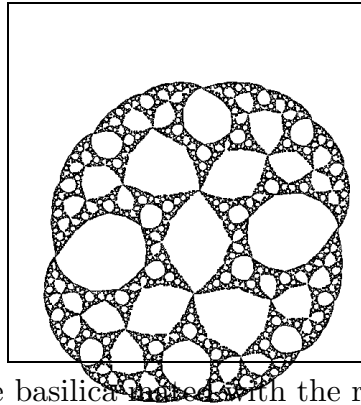

The basilica mated with the rabbit

**Question 1.** Which matings correspond to rational functions? There are some known obstructions. For example, Tan Lei has shown that matings between post-critically finite quadratic polynomials can exist only if and only if they do not belong to complex conjugate limbs of the Mandelbrot set.

**Question 2.** Can matings be constructed with quasiconformal surgery? Tan Lei uses Thurston's topological characterization of rational maps to do this. It would be nice to have a cut and paste type of construction, giving results for the case when the orbit of critical points is not finite.

**Question 3.** If one polynomial is held fixed and the other is varied continuously, does the resulting rational function vary continuously? Is mating a continuous function of two variables?

The second type of topological surgery is **tuning**. First take a polynomial $P_1$ with a periodic critical point $\omega$ of period $k$, and assume that no other critical points are in the entire basin of this superattractive cycle. Let $P_2$ be a polynomial with one critical point whose degree is the same as the degree of $\omega$. We also assume that the Julia sets of $P_1$ and $P_2$ are connected. We assume the closure $\bar{B}$ of the immediate basin of $\omega$ is homeomorphic to the closed unit disk $\bar{D}$, and that the Julia set for $P_2$ is locally connected. Now, $P_1^k$ maps $\bar{B}$ to itself by a map which is conjugate to the map $z \mapsto z^d$ of $\bar{D}$, where $d$ is the degree of the critical point. (In fact, if $d > 2$, then there are $d - 1$ possible choices for the conjugating homeomorphism, and we must choose one of them.) Intuitively the idea is now the following. Replace the basin $B$ by a copy of the dynamical plane for $P_2$, gluing the "circle at infinity" for this plane onto the boundary of $B$ so that external angles for $P_2$ correspond to internal angles in $\bar{B}$. Now shrink each external ray for $P_2$ to a point. Also, make an analogous modification at each pre-image of $B$. The map from the modified $B$ to its image will be given by $P_2$, and the map on all other inverse images of the modified $B$ will be the identity. The result,$P_3$, called $P_1$ tuned with $P_2$ at $\omega$, should be conjugate to a polynomial having the same degree as $P_1$. Conversely $P_2$ is said to be obtained from $P_3$ by renormalization.

In the case of quadratic polynomials, the tunings can be made also in the case when $P_2$ is not locally connected.

As an example we can take $P_1$ to be the rabbit polynomial. Then we can take $P_2(z) = z^2 - 2$ which has the closed segment from -2 to 2 as its Julia set. The following figure shows the resulting quadratic Julia set tuning the rabbit with the segment ($z^2 + c$ where $c \sim -.101096 + .956287i$).

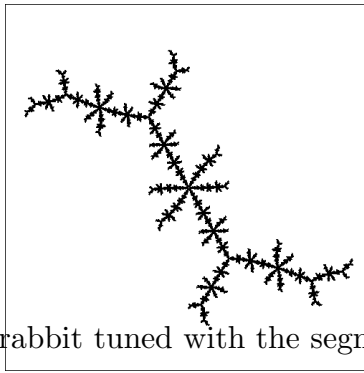

The rabbit tuned with the segment

In the picture we see each ear of the rabbit replaced with a segment.

**Question 4.** Does the tuning construction always give a result which is conjugate to a polynomial? This is true when $P_1$ and $P_2$ are quadratic.

**Question 5.** Can tunings be constructed with quasiconformal surgery?

**Question 6.** Does the resulting polynomial vary continuously with $P_2$? This is true when $P_1$ and $P_2$ are quadratic [**DH**2].

**Question 7.** Does the resulting tuning vary continuously with $P_1$? (here we consider only polynomials $P_1$ of degree greater than 2 with a superstable orbit of fixed period.)

**Question 8.** Let $P_{1,k}$ be a sequence of polynomials with a superstable orbit whose period tends to infinity. If $P_{1,k}$ tends to a limit $P_{1,\infty}$, do the tunings of $P_2$ with $P_{1,k}$ also tend to $P_{1,\infty}$?

The third kind of surgery is **intertwining surgery**.

Let $P_1$ be a monic polynomial with connected Julia set having a repelling fixed point $x_0$ which has a ray landing on it with combinatorial rotation number $p/q$. Look at the cycle of $q$ rays which are the forward images of the first. Cut along these rays and we get $q$ disjoint wedges. Now let $P_2$ be a monic polynomial with a ray of the same combinatorial rotation number landing on a repelling periodic point of some period dividing $q$ (such as 1 or $q$). Slit this dynamical plane along the same rays making holes for the wedges. Fill the holes in by the corresponding wedges above making a new sphere. The new map will be given by $P_1$ and $P_2$ except on a neighborhood of the inverse images of the cut rays where it will have to be adjusted to make it continuous. This construction should be possible to do quasiconformally using the methods in [**BD**] together with Shishikura's new (unpublished) method of presurgery in the case where the rays in the $P_2$ space land at a repelling orbit. This construction doesn't seem to work when the rays land at a parabolic orbit.

For instance we can take $P_1(z) = z^2$ and $P_2(z) = z^2 - 2$. The Julia set for $P_1$ is the unit circle with repelling fixed point at 1 and the ray at angle 0 lands on it with combinatorial rotation number 0. The Julia set for $P_2$ is the closed segment from -2 to 2 with repelling fixed point 2 and the ray at angle 0 lands on it with combinatorial rotation number 0. We cut along the 0 ray in both cases. Opening the cut in the first dynamical space gives us one wedge. The space created by opening the cut in the second space is the hole into which we put the wedge. The resulting cubic Julia set is shown in the following picture (the polynomial is $z^3 + az$ where $a \sim 2.55799i$).

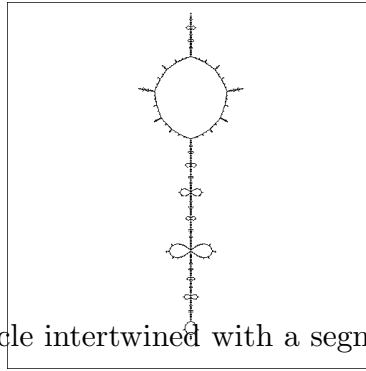

A circle intertwined with a segment

We see in the picture the circle and the segment, and at the inverse image of the fixed point on the segment we see another circle. At the other inverse of the fixed point on the circle we see a segment attached. All the other decorations come from taking various inverses of the main circle and segment.

As a second example we can intertwine the basilica with itself. The ray $1/3$ lands at a fixed point and has combinatorial rotation number $1/2$. The following is the Julia set for the basilica intertwined with itself (the polynomial here is $z^3 - \frac{3}{4}z + \frac{\sqrt{-7}}{4}$).

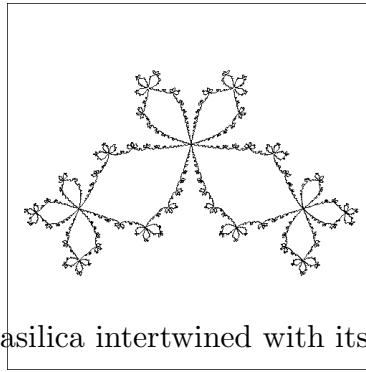

A basilica intertwined with itself

**Question 9.** When does an intertwining construction give something which is conjugate to a polynomial?

**Question 10.** Can intertwinings be constructed with quasiconformal surgery?

**Question 11.** Does the resulting polynomial vary continuously in $P_2$?

# Rational maps and Teichmüller space
## Curt McMullen

Let $X$ be a complex manifold and let $f : X \times \widehat{\mathbb{C}} \to \widehat{\mathbb{C}}$ be a holomorphic map. Then $f$ describes a family $f_\lambda(z)$ of rational maps from the Riemann sphere to itself, depending holomorphically on a complex parameter $\lambda$ ranging in $X$.

By [MSS], there is an open dense set $X_0 \subset X$ on which the family is structurally stable near the Julia set: in fact $f_a$ and $f_b$ are quasiconformally conjugate on their respective Julia sets whenever $a$ and $b$ lie in the same component $U$ of $X_0$. The mappings in $X_0$ are said to be $J$-stable.

In this note we will record some problems concerning the boundaries of components $U$, and consequently concerning limits of quasiconformal deformations of a given rational map.

**Example I. Quadratic polynomials.** The most famous such problem is the following. Let $X = \mathbb{C}$, and let $f_\lambda(z) = z^2 + \lambda$. Then $X_0$ contains a unique unbounded component $U$.

**Problem.** Is the boundary of $U$ locally connected?

This is equivalent to the question:

> Is the Mandelbrot set $M$ locally connected?

Indeed, $X_0$ is just the complement of the boundary of the Mandelbrot set.

The importance of this question is twofold. First, if $M$ is locally connected, then existing work provides detailed information about its combinatorial structure, and one has a good understanding of the "bifurcations" of a quadratic polynomial and many related maps. Secondly, the local connectivity of $M$ implies the density of hyperbolic dynamics ("Axiom A") for degree two polynomials, another well-known conjecture which has eluded proof for many years. For more details see [Dou1], [Dou2], [DH1], [DH2], [Lav], [Th].

**Compactifying the space of proper maps.** We now turn to a second example motivated by an analogy with Bers' embedding of Teichmüller space. Let $A$ and $B$ be two proper holomorphic maps of the unit disk $\Delta$ to itself, both of degree $n > 1$ and fixing zero. ($A$ and $B$ are finite Blaschke products.) Then it is well-known that $A$ and $B$ are topologically conjugate on the unit circle $S^1$, and the conjugacy $h$ is unique once we have chosen a pair of fixed points $(a, b)$ for $A$ and $B$ such that $h(a) = b$. Moreover $h$ is quasisymmetric; this is a general property of conjugacies between expanding conformal dynamical systems [Sul].

Now glue two copies of the disk together by $h$ and transport the dynamics of $A$ and $B$ to the resulting Riemann surface, which is a sphere. We obtain in this way an expanding (i.e. hyperbolic) rational map $f(A, B)$. The Julia set $J$ of $f(A, B)$ is a quasicircle, and $f$ is holomorphically conjugate to $A$ and $B$ on the components of the complement of $J$. The mapping $f(A, B)$ is determined by $h$ up to conformal conjugacy.

We will loosely speak of spaces of mappings as being "the same" if they represent the same conformal conjugacy classes. It is often useful to require that the conjugacy preserves some finite amount of combinatorial data, such as a distinguished fixed point. For simplicity we will gloss over such considerations below.

**Example II.** Let $X$ be the space of degree $n$ polynomials, $X_0$ the open dense subset of $J$-stable polynomials and $U$ the component of $X_0$ containing $z^n$. Then $U$ is the same as the set of maps of the form $f(z^n, B)$. Equivalently, $U$ consists of those polynomials with an attracting fixed point with all critical points in its immediate basin.

Let us denote this set of polynomials by $\mathcal{B}(z^n)$. It is easy to see that $\mathcal{B}(z^n)$ is an open set of polynomials with compact closure. Thus this construction supplies both a complex structure for the space of Blaschke products, and a geometric compactification of that space.

**Problem.** Describe the boundary of $\mathcal{B}(z^n)$ in the space of polynomials of degree $n$.

For degree $n = 2$ this is easy (the boundary is a circle) but for $n = 3$ it is already subtle.

To explain the kind of answer one might expect, we consider not one boundary but many. More precisely, let $\mathcal{B}(A)$ denote the space of rational maps $f(A, B)$ for some other *fixed* $A$ and varying $B$. This space also inherits a complex structure and the map $f(z^n, B) \mapsto f(A, B)$ gives an biholomorphic map

$$F : \mathcal{B}(z^n) \to \mathcal{B}(A).$$

The closure of $\mathcal{B}(A)$ is the space of rational maps provides another geometric compactification of this complex manifold.

**Problem.** Show that for $n > 2$ and $A \neq z^n$, $F$ does *not* extend to a homeomorphism between the boundaries of $\mathcal{B}(z^n)$ and $\mathcal{B}(A)$.

Thus we expect that the complex space $\mathcal{B}$ (whose complex structure is independent of $A$) has many natural geometric boundaries. But perhaps the lack of uniqueness can be accounted for by the presence of *complex submanifolds* of the boundary, i.e.

by the presence of rational maps in the compactification which admit quasiconformal deformations.

To make this precise, let $\partial(A)$ denote the quotient of the boundary of $\mathcal{B}(A)$ by the equivalence relation $f \sim g$ if $f$ and $g$ are quasiconformally conjugate (equivalently, if $f$ and $g$ lie in a connected complex submanifold of the boundary). The resulting space (in the quotient topology) still forms a boundary for $\mathcal{B}(A)$, but it is non-Hausdorff when $n > 2$.

**Conjecture.** The holomorphic isomorphism $F : \mathcal{B}(z^n) \to \mathcal{B}(A)$ extends to a homeomorphism from $\partial(z^n)$ to $\partial(A)$.

**Problem.** Give a combinatorial description of the topological space $\partial(z^n)$.

Such a description may involve laminations, as discussed in [Th].

**An analogy with Teichmüller theory.** The "mating" of $A$ and $B$ has many similarities with the mating of Fuchsian groups uniformizing a pair of compact genus $g$ Riemann surfaces $X$ and $Y$. Such a mating is provided by Bers' simultaneous uniformization theorem [Bers]. The result is a Kleinian group $\Gamma(X, Y)$ whose limit set is a quasicircle. Moreover, fixing $X$, the map $Y \mapsto \Gamma(X, Y)$ provides a holomorphic embedding of the Teichmüller space of genus $g$ into the space of Kleinian groups. One can then form a boundary for Teichmüller space by taking the closure.

It has recently been shown that *this* boundary does indeed *depend* on the base point $X$ [KT]. However Thurston has conjectured that the space $\partial(X)$, obtained by identifying quasiconformally conjugate groups on the boundary, is a (non-Hausdorff) boundary which is independent of $X$.

Moreover a combinatorial model for $\partial(X)$ is conjecturally constructed as follows. Let $\mathbb{P}\mathcal{ML}$ denote the space of projective measured laminations on a surface of genus $g$; then $\partial(X)$ is homeomorphic to the quotient of $\mathbb{P}\mathcal{ML}$ by the equivalence relation which forgets the measure. (See [FLP] for a discussion of $\mathbb{P}\mathcal{ML}$ as a boundary for Teichmüller space.)

**Remarks.**

1. We do not expect that one can give a combinatorial description of the "actual" boundary of $\mathcal{B}(z^n)$ (in the space of polynomials). For similar reasons, we believe it unlikely that one can describe the uniform structure induced on the space of critically finite rational maps by inclusion into the space of all rational maps.

2. It is known that Teichmüller space is a domain of holomorphy. So it is natural to ask the following intrinsic:

**Question.** Is $\mathcal{B}(z^n)$ a domain of holomorphy? More generally, is every component of the space of expanding rational maps (or polynomials) a domain of holomorphy?

**Density of cusps.** The preceding discussion becomes interesting only when the space of rational maps under consideration has two or more (complex) dimensions. We conclude with two concrete questions about boundaries in a one-parameter family of rational maps.

**Example III.** Let

$$f_\lambda(z) = \lambda z^2 + z^3$$

where $\lambda$ ranges in $X = \mathbb{C}$, and let $U$ denote the component of $X_0$ containing the origin. That is, $U$ is the set of $\lambda$ for which both finite critical points are in the immediate basin of zero.

A *cusp* on $\partial U$ is an $f_\lambda$ with a parabolic periodic cycle.

**Conjecture.** Cusps are dense in $\partial U$.

This conjecture is motivated by the density of cusps on the boundary of Teichmüller space [Mc]. It is not hard to show that it is implied by the following:

**Conjecture.** The boundary of $U$ is a Jordan curve.

# Thurston's algorithm without critical finiteness
## John Milnor

Thurston's algorithm is a powerful method for passing from a topological branched covering $S^2 \to S^2$ to a rational map having closely related dynamical properties. (See [DH].) The same method can be used to pass from a piecewise monotone map of the interval to a closely related polynomial map of the interval.

Suppose that we start with an orientation preserving branched covering map $f_0 : S^2 \to S^2$. We identify $S^2$ with the Riemann sphere $\overline{\mathbf{C}} = \mathbf{C} \cup \infty$. In order to anchor this sphere, choose three base points. (For best results, choose dynamically significant base points, for example periodic points of $f_0$, or critical points, or critical values.)

> **Lemma:** *There is one and only one homeomorphism $h_0 : S^2 \to S^2$ which fixes the three base points, and which has the property that the composition $r_0 = f_0 \circ h_0$ is holomorphic, or in other words is a rational map.*

[Proof: Let $\sigma_0$ be the standard conformal structure on the 2-sphere, and let $\sigma = f_0^*(\sigma_0)$ be the pulled back conformal structure, so that $f_0$ maps $(S^2, \sigma)$ holomorphically onto $(S^2, \sigma_0)$. Then $h_0$ must be the unique conformal isomorphism from $(S^2, \sigma_0)$ onto $(S^2, \sigma)$ which fixes the three base points.] Now consider the map $f_1 = h_0^{-1} \circ f_0 \circ h_0$, which is topologically conjugate to $f_0$. In this way, we obtain a commutative diagram

$$
\begin{array}{ccc}
S^2 & \xrightarrow{\ f_1\ } & S^2 \\[4pt]
h_0 \downarrow & r_0 \searrow & h_0 \downarrow \\[4pt]
S^2 & \xrightarrow{\ f_0\ } & S^2 \ .
\end{array}
$$

Continuing inductively, we produce a sequence of branched coverings $f_n$, and a sequence of homeomorphisms $h_n$ fixing the base points, so that $f_{n+1} = h_n^{-1} \circ f_n \circ h_n$, and so that each composition $r_n = f_n \circ h_n$ is a rational map. The marvelous property of this construction is that in many cases the homeomorphisms $h_n$ seem to tend uniformly to the identity, so that the successive maps $f_n$, which are all topologically conjugate to $f_0$, come closer and closer to the rational maps $r_n$. In fact the sequence of compositions $\phi_n = (h_0 \circ \cdots \circ h_n)^{-1}$ may converge uniformly to a limit map $\phi$, at least on the non-wandering set. In this case, it follows that the rational limit map is topologically semi-conjugate (or perhaps even conjugate) to $f_0$,

$$
r_\infty \circ \phi \ = \ \phi \circ f_0
$$

on the non-wandering set.

> **Problem:** *Under what conditions will this sequence of rational maps $r_n$ converge uniformly to a limit map $r_\infty$ ? Under what conditions, and on what subset of $S^2$, will the maps $\phi_n$ converge uniformly to a limit?*

In the post-critically finite case, Thurston defines an obstruction, which vanishes if and only if the restriction of the $\phi_n$ to the post-critical set converges uniformly to a one-to-one limit function. If this obstruction vanishes, then it follows that the $r_n$ converge.

However, there would be interesting applications where $f_0$ is not post-critically finite, so that no such criterion is known. A typical example is provided by the problem of "mating". (Compare Bielefeld's discussion, as well as [Ta], [Sh].) Let $p$ and $q$ be monic polynomial maps having the same degree $d \geq 2$. Conjugating $p$ by the diffeomorphism $z \mapsto z/\sqrt{1 + |z|^2}$ from $\mathbf{C}$ onto the unit disk $D$, we obtain a map $p^*$ which extends smoothly over the closed disk $\bar{D}$. Similarly, conjugating $q$ by $z \mapsto \sqrt{1 + |z|^2}/z$ we obtain a map $q^*$ which extends smoothly over the complementary disk $\mathbf{C} \smallsetminus D$. Now $p^*$ and $q^*$ together yield a $C^1$-smooth map $f_0 : \bar{\mathbf{C}} \to \bar{\mathbf{C}}$, and we can apply Thurston's method as described above. If this procedure converges to a well behaved limit, then the resulting rational map $r_\infty$ of degree $d$ may be called the "mating" of $p$ and $q$.

**Maps of the interval.** The situation here is quite similar. Let $f_0$ be a piecewise-monotone map of the interval $I = [0, 1]$ with $d$ alternately ascending and descending laps, and suppose that $f_0$ carries the boundary points 0 and 1 to boundary points. Then there is one and only one orientation preserving homeomorphism $h_0$ of the interval such that the composition $p_0 = f_0 \circ h_0$ is a polynomial map of degree $d$. Setting $f_1 = h_0^{-1} \circ f_0 \circ h_0$, we can proceed inductively, constructing homeomorphisms $h_n$, polynomials $p_n = f_n \circ h_n$, and topologically conjugate maps $f_{n+1} = h_n^{-1} \circ f_n \circ h_n$. Again the problem is to decide when and where this procedure converges.

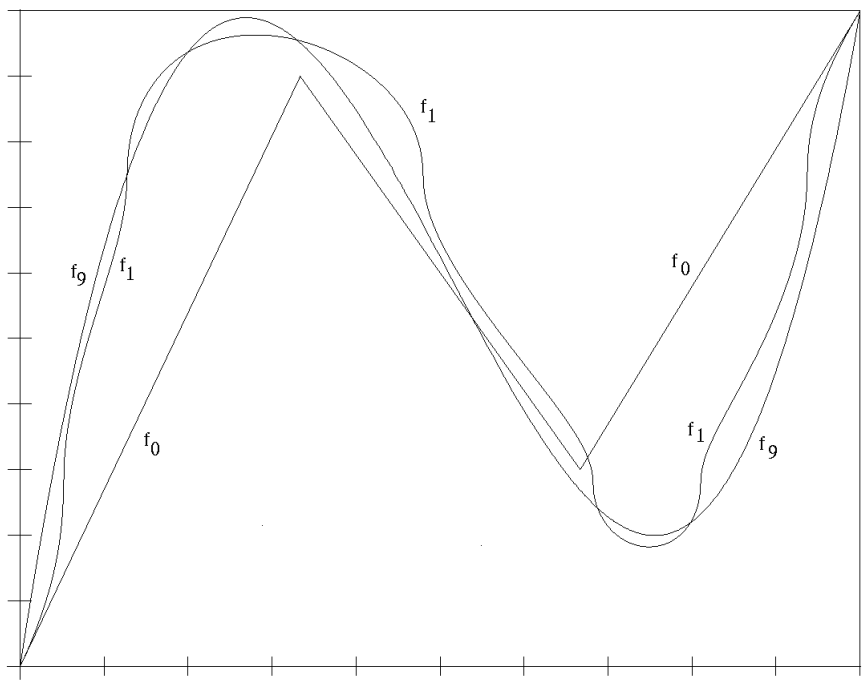

*A typical run of Thurston's method, starting with a piece-wise linear map $f_0$ of the interval. (Horizontal scale exaggerated.) The graphs of $f_0$, $f_1$ and $f_9$ are shown. The latter seems indistinguishable from $f_\infty = p_\infty$.*

**References.**

Stony Brook, April 1992

# A Possible Approach to a Complex Renormalisation Problem
## Mary Rees

**Preliminary Definitions**. For a branched covering $f : \overline{\mathbf{C}} \to \overline{\mathbf{C}}$, we define

$$X(f) = \{f^n(c) : c \text{ critical}, \ n > 0\}.$$

Then $f$ is **critically finite** if $\#(X(f))$ is finite. Two critically finite branched coverings $f_0$, $f_1$ are **(Thurston) equivalent** if there is a path $f_t$ through critically finite branched coverings connecting them with $X(f_t)$ constant in $t$.

We are only concerned, here, with orientation-preserving degree two branched coverings for which one critical point is fixed and the other is periodic. By a theorem of Thurston's ([**T**], [**D-H**]), any such branched covering $f_0$ is equivalent to a unique degree two polynomial $f_1$ of the form $z \mapsto z^2 + c$ (some $c \in \mathbf{C}$).

Now let $f_1$, $f_2$ be two degree two polynomials of the form $z \mapsto z^2 + c_i$ ($i = 1$, $2$), with $0$ periodic of periods $m$, $n$ respectively. Then we define **the tuning of $f_1$ about $0$ by $f_2$**, written $f_1 \vdash f_2$, as follows. This is simply a branched covering defined up to equivalence. Let $D$ be an open topological disc about $0$ such that the discs $f_1^i(D)$ ($0 \leq i < m$) are all disjoint, $f_1^m(D) \subset D$ and $f_1 : f_1^i(D) \to f_1^{i+1}(D)$ is a homeomorphism for $1 \leq i < m$. Let $g$ be a rescaling of $f_2$, and $V$ a closed bounded topological disc with $V \subset gV \subset f_1^m(D)$ whose complement is in the attracting basin of $\infty$ for $g$. Then we define

$$f_1 \vdash f_2 \quad \begin{aligned} &= f_1 \text{ outside } D, \\ &= f_1^{-(m-1)} \circ g \text{ in } V, \end{aligned}$$

and extend to map the annulus $D \setminus V$ by a two-fold covering to $f_1^m(D) \setminus g(V)$. Then

$$(f_1 \vdash f_2)^m = g \text{ in } V.$$

Thus $f_1 \vdash f_2$ is critically finite with $0$ of period $n \cdot m$, and is equivalent to a unique polynomial $z \mapsto z^2 + c$.

For any sequence $\{f_i\}$ of polynomials, we can also define $f_1 \vdash \cdots \vdash f_n$ for all $n$.

For concreteness, we consider the following renormalisation problem, but different versions are possible.

Let $\{f_i\}$ be any sequence of polynomials of the form $z \mapsto z^2 + c_i$, where the $f_i$ (and $c_i$) take only finitely many different values, and $0$ is of period $m_i$ under $f_i$. Write $g_n$ for the polynomial $z \mapsto z^2 + c$ equivalent to $f_1 \vdash \cdots \vdash f_n$. and

$$n_k = \prod_{i \leq k} m_i.$$

**Problem.** Prove geometric properties of $X(g_n)$. Specifically, show that the set

$$\{g_n^{n_k \ell + i}(0) : 0 \leq \ell < m_{k+1}\} \tag{1}$$

has uniformly bounded geometry for all $i \leq n_k$, $k < n$ and all $n$.

Of course, this problem (and stronger versions) is not new, has been the focus of much effort, and, in the real case, has been resolved by Sullivan [**S**]. The most obvious method of approach (which was not, in the end, efficacious in the real case) is through analysis of the main technique used to prove Thurston's theorem mentioned above. We now recall this.

**Thurston's Pullback Map on Teichmüller space.**

To simplify, we stick to orientation-preserving degree two critically finite branched coverings with fixed critical value $v_2$ and periodic critical value $v_1$. Let $g$ be one such. Let $X = X(g)$. We let $s : \overline{\mathbf{C}} \to \overline{\mathbf{C}}$ be given by $s(z) = z^2$. Let $\mathcal{T} = \mathcal{T}(X)$ be the Teichmüller space of the sphere with set of marked points $X$, so that

$$\mathcal{T} = \{[\varphi] : \varphi \text{ is a homeomorphism of } \overline{\mathbf{C}}\}$$

and $[\varphi]$ denotes the quotient of the isotopy class under isotopies constant on $X$ by left Möbius composition, that is, $[\varphi] = [\sigma \circ \varphi \circ \psi]$ for any Möbius transformation $\sigma$ and $\psi$ isotopic to the identity rel $X$. It is convenient to choose representatives $\varphi$ so that $\varphi(v_1) = 0$, $\varphi(v_2) = \infty$. Then

$$\tau : \mathcal{T} \to \mathcal{T}$$

is defined by

$$\tau([\varphi]) = [s^{-1} \circ \varphi \circ g].$$

(The righthand side makes perfectly good sense as a homeomorphism.)

By Thurston's theorem (in this setting), $\tau$ is a contraction with respect to the Teichmüller metric $d$ on $\mathcal{T}$, and has a unique fixed point $[\varphi]$. Then there are $\psi$ isotopic to the identity via an isotopy fixing $X(g)$ (unique given $\varphi$) and a Möbius transformation $\sigma$ such that

$$g = \varphi^{-1} \circ s \circ \sigma \circ \varphi \circ \psi.$$

In particular, $g$ and $s \circ \sigma$ are equivalent, and $X(s \circ \sigma) = \varphi(X(g))$.

**The "Obvious" Method of Approach.**

We can choose $h_n$ equivalent to $g_n$ so that the sets of (1), with $h_n$ replacing $g_n$, have uniformly bounded geometry for $i \leq n_k$, $k < n$, and all $n$. Then let $\mathcal{T}_n = \mathcal{T}(X(h_n))$, and $\tau_n : \mathcal{T}_n \to \mathcal{T}_n$ be the associated pullback. It suffices (!) to prove convergence, as $m \to \infty$, and uniform in $n$, of the sequences $\{\tau_n^m(\text{identity})\}$. Of course, for fixed $n$, the convergence would be with respect to the Teichmüller metric $d_n$ on $\mathcal{T}_n$. This seems to be impossible to implement. An alternative is suggested below. One virtue - and probably the only one - of this alternative is that it has not yet been tried (so far as I know). Before making this precise, we need to clarify some properties of the Teichmüller metric.

**The Teichmüller metric and its Derivative.**

Let $\mathcal{T} = \mathcal{T}(X)$ (for any finite set $X \subset \overline{\mathbf{C}}$) and let $d$ denote the Teichmüller metric. Let $[\varphi]$, $[\psi] \in \mathcal{T}$. Assume without loss of generality that $\infty \in \varphi(X)$, $\psi(X)$. Then there is a unique quasiconformal homeomorphism $\chi : \overline{\mathbf{C}} \to \overline{\mathbf{C}}$ with the following properties.

1. $\chi(\varphi(X)) = \psi(X)$ and $[\chi \circ \varphi] = [\psi]$.

2. There is a rational function $q$ with at most simple poles in $\mathbf{C}$, all occurring at points of $\varphi(X)$, and at least three more poles than zeros in $\mathbf{C}$, such that the directions of maximal stretch and contraction of $\chi$ are tangent to the vector fields $i\sqrt{\overline{q}}$, $\sqrt{\overline{q}}$ respectively, and the dilatation (ratio of infinitesimal stretch to contraction) is constant.

3. The images under $\chi_*$ of these vector fields are of the form $i\sqrt{\overline{p}}$, $\sqrt{\overline{p}}$, for a rational function $p$ with at most simple poles in $\mathbf{C}$, all occurring at points of $\psi(X)$.

The function $q$ is then also unique, up to a positive scalar multiple, and becomes unique if we normalise so that

$$\int \mid q \mid \frac{d\overline{z} \wedge dz}{2i} = 1.$$

Similarly, we normalise $p$. (Of course, $q$ represents a quadratic differential $q(z)dz^2$, but it is convenient to keep the representing rational function in the foreground.)

Let $h = (h(x)) \in \mathbf{C}^X$ be small, taking $h(x) = 0$ if $\varphi(x) = \infty$. Then by abuse of notation, we write $\varphi + h$ for a homeomorphism near $\varphi$ with $(\varphi + h)(x) = \varphi(x) + h(x)$. Then the following holds, where $q$, $p$ are demined by $[\varphi]$, $[\psi]$ as above $[\mathbf{R}]$.

$$d([\varphi + h], [\psi + k]) = d([\varphi], [\psi]) + 2\pi \operatorname{Re}\left(\sum_{x \in X}\left(\operatorname{Res}(q, \varphi(x))h(x) - \operatorname{Res}(p, \psi(x))k(x)\right)\right)$$
$$+ o(h) + o(k).$$

Now we consider the case $X = X(g)$ and $y = \tau x$. As before we consider only specific $g$ and take $s(z) = z^2$ (as before). The **pushforward** $s_*q$ of a rational function (or quadratic differential) $q$ is defined by

$$s_*q(z) = \sum_{s(w) = z} \frac{q(w)}{s'(w)^2}$$

if $s'(w) \neq 0$ for $s(w) = z$. If $q$ has only simple poles in $\mathbf{C}$, and at least 3 more poles than zeros in $\mathbf{C}$, then $s_*q$ extends to a rational function on $\mathbf{C}$ with the same properties. Then if $q$, $p$ are determined by $[\varphi]$, $\tau([\varphi])$, taking $\varphi(v_1) = 0$, $\varphi(v_2) = \infty$ as above,

$$d([\varphi + h], \tau([\varphi + h])) = d([\varphi], \tau([\varphi])) + 2\pi \operatorname{Re}\left(\sum_{x \in X}\left(\operatorname{Res}(q, \varphi(x)) - \operatorname{Res}(s_*p, \varphi(x))\right)h(x)\right)$$
$$+ o(h).$$

**The Suggested Alternative Approach to the Problem.** Take $\mathcal{T} = \mathcal{T}(X(g))$, $\tau : \mathcal{T} \to \mathcal{T}$. Let

$$F([\varphi]) = d([\varphi], \tau([\varphi])).$$

Then the derivative formula for $F$ above theoretically enables us to construct flows for which $F$ decreases along orbits. It can be shown that the only critical point of $F$ occurs where $F = 0$. So if we can find a compact subset $B$ of $\mathcal{T}$ with smooth boundary and a vector field $v$ pointing inward on $\partial B$ with $DF(v) < 0$, then the (unique) zero of $F$ must be inside $B$. Now put a subscript $n$ on everything. Conceivably we can find $B_n \subset \mathcal{T}_n$ and vector field $v_n$ pointing inward on $\partial B_n$ such that if $A \subset X(g_n)$ is any of the sets in (1) and $[\varphi] \in \partial B_n$ then $\varphi(A)$ has bounded geometry (uniformly in $A$, $n$) and $DF_n(v_n) < 0$?

# Section 2: Geometry of Julia Sets

## Geometry of Julia sets

### Lennart Carleson

The geometry of connected Julia sets for hyperbolic quadratic polynomials is now well understood. Bounded components of the Fatou set are quasi-circles while the unbounded component is a John domain.

The geometry of a flower (for a rational fixed point) is also known. If the flower has more than one petal, each component is a quasi-disk. The 1-petal flower is a John domain (see a forthcoming paper by P.Jones and L.Carleson in Boletin de Brasil).

For Siegel disks $S$, a basic result by M. Herman is that the critical point belongs to the boundary of $S$ if the rotation number $\lambda = e^{2\pi i\theta}$ is a Siegel number, i.e.

$$|\theta - \frac{p}{q}| > \frac{c}{q^n} \text{ for some } c > 0, \ n < \infty,$$

or more generally when the arithmetic condition of J.-C. Yoccoz's global theorem on conjugacy of analytic diffeomorphisms of the circle is satisfied.

Another remarkable result of M.Herman is that when the critical point belongs to $\partial S$, then $S$ is a quasi-disk if and only if $\lambda$ is of bounded type, i.e.

$$|\theta - \frac{p}{q}| > \frac{c}{q^2}.$$

With J.-C. Yoccoz he has also proved that $\partial S$ is a Jordan curve for almost all $\theta$. It is not known which arithmetic condition implies this. E.g., is there $\gamma_0 > 2$ so that $|\theta - p/q| > C/q^\gamma$ implies that $S$ is a Jordan domain for $\gamma < \gamma_0$ but not for $\gamma > \gamma_0$?

A particularly interesting question concerns the geometry of $\partial S$ at the critical point. Computer experiments show that in many cases $\partial S$ has an angle of about $120°$ opening at the critical point. Prove this at least for $\theta = \theta_0 = (\sqrt{5} - 1)/2$. For this value there should also exist a renormalization at the critical point.

There is also a very interesting regularity of the Taylor coefficients of the conjugating map. Consider more generally the family $\mathcal{P}_\rho(z)$, $\mathrm{Re}(\rho) > 0$, with

$$\mathcal{P}'_\rho = \lambda(1-z)^\rho, \quad \mathcal{P}_\rho(0) = 0$$

so that $\rho = 1$ corresponds to $\lambda(z - z^2/2)$. Let $h(\zeta)$ be the conjugating map in $|\zeta| < 1$ with $h(1) = 1$. (For general $\rho$ the proof that $1 \in \partial S$ is not known, but should be rather similar to the case $\rho = 1$). Form

$$f(\zeta) = \frac{h'(\zeta)}{1 - h(\zeta)} = \sum_0^\infty a_\nu \zeta^\nu.$$

Then

$$f' - f^2 = f\rho \sum_0^\infty (\frac{1}{2} + \frac{i}{2}\cot(\nu+1)\pi\theta)a_\nu\zeta^\nu.$$

If the imaginary part in the parenthesis is dropped we obtain

$$f'_0 = (1 + \frac{\rho}{2})f_0^2, \quad f_0 = \frac{1}{(1+\rho/2)(1-z)}, \quad h_0 = (1-z)^{2/(\rho+2)}.$$

Computer experiments indicate for $\theta = \theta_0$, $\rho = 1$

$$|a_\nu - \frac{2}{3}| < 0.1 \quad \text{(say) for all } \nu,$$

where $2/3$ corresponds to $f_0$. It would be interesting to make the approximation rigorous at least for small $\rho$.

In the non-hyperbolic case very little is known (and very little can be probably said in general). The simplest case of a strictly preperiodic critical point leads to John domains (the Julia set is called a dendrite). It should be possible to analyse the general Misiurewicz case when the critical point never returns close to itself. In the case of $1 - az^2$, $a$ is real, this condition is equivalent to the Fatou set being a John domain. To which extent does this hold for general Misiurewicz points?

# Problems on local connectivity.[1]

## John Milnor

If the Julia set $J(f)$ of a quadratic polynomial is connected, then Yoccoz has proved[2] that $J(f)$ is locally connected, unless either:

(1) $f$ has an irrationally indifferent periodic point, or

(2) $f$ is infinitely renormalizable.

**Cremer Points.** To illustrate case (1), consider the polynomial

$$P_\alpha(z) = z^2 + e^{2\pi i\alpha}z$$

with a fixed point of multiplier $\lambda = e^{2\pi i\alpha}$ at the origin. Take $\alpha$ to be real and irrational. For generic choice of $\alpha$ (in the sense of Baire category), Cremer showed that there is no local linearizing coordinate near the origin. We will say briefly that the origin is a *Cremer point*, or that $P_\alpha$ is a *Cremer polynomial*. According to Sullivan and Douady, the existence of such a Cremer point implies that the Julia set is not locally connected. More explicitly, let $t(\alpha)$ be the angle of the unique external ray which lands at the corresponding point of the Mandelbrot set. For generic choice of $\alpha$, Douady has shown that the corresponding ray in the dynamic plane does not land, but rather has an entire continuum of limit points in the Julia set. (Compare [Sø].) Furthermore, the $t(\alpha)/2$ ray in the dynamic plane accumulates both at the fixed point 0 and its pre-image $-\lambda$.

**Problem 1.** Is there an arc joining 0 to $-\lambda$, in the Julia set of such a Cremer polynomial?

**Problem 2.** Give a plausible topological model for the Julia set of a Cremer polynomial.

**Problem 3.** Make a good computer picture of the Julia set of some Cremer polynomial.

**Problem 4.** Can there be any external rays landing at a Cremer point?

**Problem 5.** Can the critical point of a Cremer polynomial be accessible from $\mathbf{C} \smallsetminus J$?

**Problem 6.** If we remove the fixed point from the Julia set of a Cremer polynomial, how many connected components are there in the resulting set $J(P_\alpha)\smallsetminus\{0\}$, ie., is the number of components countably infinite?

**Problem 7.** The Julia set for a generic Cremer polynomial has Hausdorff dimension two. Is this true for an arbitrary Cremer polynomial? Do Cremer Julia sets have measure zero? (Compare [Sh], [L1], [L2].)

In the quadratic polynomial case, Yoccoz has shown that every neighborhood of a Cremer point contains infinitely many periodic orbits. On the other hand, Perez-Marco [P-M1] has described non-linearizable local holomorphic maps for which this is not true.

**Problem 8.** For a Cremer point of an arbitrary rational map, does every neighborhood contain infinitely many periodic orbits?

---

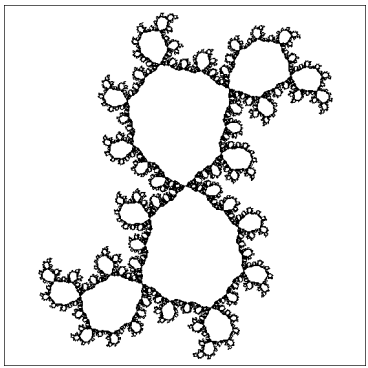

*Figure 1. Julia set of $P_\alpha$ where $\alpha = .78705954039469$ has been randomly chosen.*

**Siegel Disks.** (Compare Carleson's discussion.) If $\alpha$ satisfies a Diophantine condition (in particular, for Lebesgue almost every $\alpha$), Siegel showed that there is a local linearizing coordinate for the polynomial $P_\alpha(z) = z^2 + e^{2\pi i \alpha} z$ in some neighborhood of the origin. Briefly we say that the origin is the center of a *Siegel disk* $\Delta$, or that $P_\alpha$ is a *Siegel polynomial*. Yoccoz has given a precise characterization of which irrational angles yield Siegel polynomials and which yield Cremer polynomials. ([Y], [P-M2].)

Herman, making use of ideas of Ghys, showed that there exists a value $\alpha_0$ so that $P_{\alpha_0}$ has a Siegel disk whose boundary $\partial \Delta$ does not contain the critical point. It follows that the Julia set $J(P_{\alpha_0})$ is not locally connected. On the other hand if $\alpha$ satisfies a Diophantine condition, then Herman showed that $\partial \Delta$ does contain the critical point.

**Problem 9.** Give any example of a Siegel polynomial whose Julia set is provably locally connected. Is $J(P_\alpha)$ locally connected for Lebesgue almost every choice of $\alpha$? (Compare Figure 1.) What can be said about the Hausdorff dimension of $J(P_\alpha)$?

**Problem 10.** Can a Siegel disk have a boundary which is not a Jordan curve?

**Problem 11.** Does any rational function have a Siegel disk with a periodic point in its boundary? Such an example would be extremely pathological. (In the polynomial case, Poirier has pointed out that at least there cannot be a Cremer point in the boundary of a Siegel disk. See [GM].)

**Infinitely Renormalizable Polynomials.** A quadratic polynomial $f_c(z) = z^2 + c$ is *renormalizable* if there exists an integer $p \geq 2$ and a neighborhood $U$ of the critical point zero so that the orbit of zero under $f^{\circ p}$ remains in this neighborhood forever, and so that the map $f^{\circ p}$ restricted to $U$ is polynomial-like of degree $2$. (Thus the closure $\bar{U}$ must contain no other critical points of $f^{\circ p}$, and must be contained in the interior of $f^{\circ p}(U)$.) Let $M$ be the Mandelbrot set, and let $H \subset M$ be any hyperbolic component of period $p \geq 2$. Douady and Hubbard [DH2] show that $H$ is contained in a small copy of $M$. This small copy is the image of a homeomorphic embedding of $M$ into itself, which I will denote by $c \mapsto H * c$. The elements of these various small copies $H * M \subset M$ (possibly with the root point $H * \frac{1}{4}$ removed) are precisely the renormalizable elements of $M$.

Now consider an infinite sequence of hyperbolic components $H_1, H_2, \ldots \subset M$. If the $H_i$ converge to the root point $1/4$ sufficiently rapidly, then Douady and Hubbard (unpublished) show that the intersection $\bigcap_k H_1 * \cdots * H_k * M$ consists of a single point $c_\infty$ such that the corresponding Julia set $J(f_{c_\infty})$ is not locally connected.

**Problem 12.** Suppose that $f_{c_\infty}$ is infinitely renormalizable of bounded type. For example, suppose that $c_\infty \in \bigcap_k H_1 * \cdots * H_k * M$, where the $H_i$ are all equal. Does it then follow that $J(f_{c_\infty})$ is locally connected? As the simplest special case, if we take $H_1 = H_2 = \cdots$ to be the period two component centered at $-1$, then $f_{c_\infty}$ will be the quadratic *Feigenbaum map*. Is the Julia set for the Feigenbaum map locally connected?

**Problem 13.** More generally, if $c$ is real (belonging to the intersection $M \cap \mathbf{R} = [-2, \, 1/4]$), does it follow that the Julia set $J(f_c)$ is locally connected?

**The Mandelbrot Set.** Here the most basic remaining question is the following.

**Problem 14.** Does every infinite intersection of the form $\bigcap_k H_1 * \cdots H_k * M$ reduce to a single point? Equivalently, is the set of infinitely renormalizable points totally disconnected? Does this set have measure zero? Does it in fact have small Hausdorff dimension?

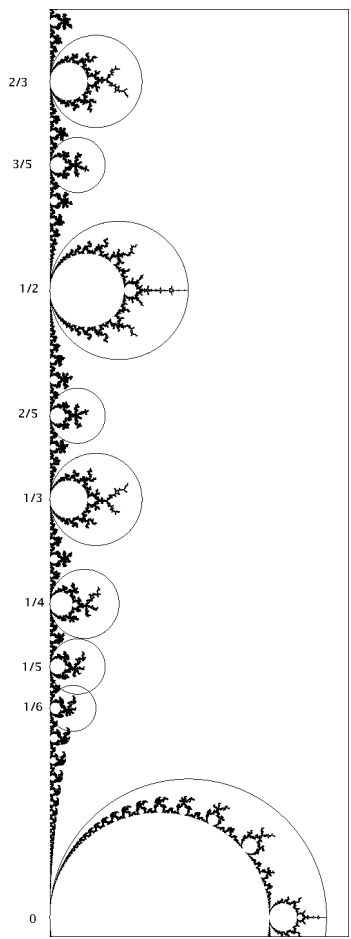

*Figure 2. Picture of the* $\log \lambda$-plane, showing the Yoccoz *disks of radius* $\log(2)/q$. (*Heights in units of* $2\pi$.)

**Problem 15.** For each rational number $0 < p/q < 1$ let $M(p/q)$ be the limb of the Mandelbrot set with interior angle $p/q$. Is the diameter of $M(p/q)$ less than $k/q^2$ for some constant $k$ independent of $p$ and $q$? If not, is it at least less than $k\log(q)/q^2$? (It is actually more natural to work in the $\log \lambda$ plane, where $f(z) = z^2 + \lambda z$. The Yoccoz inequality asserts that the corresponding limb in this $\log \lambda$ plane is contained in a disk of radius $\log(2)/q$. Compare [P], and see Figure 2.)

# Section 3: Measurable Dynamics

## Measure and Dimension of Julia Sets
### Mikhail Lyubich

**Problem 1.** Can it happen that a nowhere dense Julia set has positive Lebesgue measure?

The corresponding Ahlfors problem in Kleinian groups is also still unsolved.

So far it is known that the Julia set has zero measure in the following cases:

(i) hyperbolic, subhyperbolic and parabolic cases [DH], [L1].

(ii) a cubic polynomial with one simple non-escaping critical point and with a "non-periodic tableaue" (McMullen, see [BH]);

(iii) a quadratic polynomial which is only finitely renormalizable and has no neutral irrational cycles (Lyubich [L2] and Shishikura (unpublished)).

Let us say that a polynomial with one non-escaping critical point $c$ is *renormalizable* if there is a quadratic-like map $f^n : U \to V$, $c \in U \subset V$ $n > 1$, with connected Julia set. It corresponds to the case of periodic tableaue. Cases (i) and (ii) can be generalized in the following way:

(iv) a polynomial of any degree but with only one non-escaping critical point which does not have irrational neutral points and which is only finitely renormalizable.

In higher degrees one can describe a wide class of combinatorics for which the Julia set has zero measure (non-recurrent and "reluctantly recurrent" cases). The basic examples for which the answer is still unclear are

1. The Feigenbaum quadratic polynomial.

2. The Fibonacci polynomial $z \mapsto z^d + c$ with $d > 2$ (see [BH] or [LM] for the definition of the Fibonacci polynomial).

3. A polynomial with a Cremer point or Siegel disk (see the disciussion in Milnor's notes).

In the case when the Julia set coincides with the whole sphere the corresponding question is the following.

**Problem 2.** Is it true for all $f$ with $J(f) = \bar{\mathbf{C}}$ that the following hold?

(i) $\omega(z) = \bar{\mathbf{C}}$ for almost all $z \in \bar{\mathbf{C}}$?

(ii) $f$ is conservative with respect to the Lebesgue measure? (Conservativity means that the Poincaré Return Theorem holds).

Note that for the interval maps (replacing $\bar{\mathbf{C}}$ by an interval on which $f$ is topologically mixing) (i) and (ii) are equivalent [BL2]. Moreovere, both of them hold for the quadratic-like maps of the interval [L3].

**Problem 3.** Let again $J(f) = \bar{\mathbf{C}}$. Is it true that $f$ is ergodic with respect to the Lebesgue measure? Is it at least true that it has at most $2 \deg f - 2$ ergodic components?

The answer to the first question is yes for a large set of rational maps [R]. The answer to the second one is yes for interval maps [BL1].

The discussed problems are closely related to the deformation theory of rational maps. The link between them is given by the notion of *measurable invariant line field on the Julia set* (see [MSS]). Each such field generates a quasi-conformal deformation of $f$ supported on the Julia set. There is a series of Lattes examples having an invariant line field on the Julia set, and in these examples $J(f) = \bar{\mathbf{C}}$. Such a phenomenon is impossible at all for finitely generated Kleinian groups [S].

**Problem 4.** (Sullivan) Are the Lattes examples the only ones having measurable invariant line fields on the Julia sets?

Let us consider now an analytic family $A$ of rational maps, and denote by $Q \subset A$ the set of $J$-unstable maps.

A recent remarkable result by Shishikura [Sh] says that in the quadratic family $z \mapsto z^2 + c$ there are a lot of Julia sets with Hausdorff dimension 2.

**Problem 5** Find an explicit example of a Julia set of Hausdorff dimension 2. What is a natural geometric measure in the case when $J(f)$ has Hausdorff dimension 2 but zero Lebesgue measure?

A more general program is to develop an appropriate Thermodynamical Formalism in non-hyperbolic situations.

**Problem 6.** (i) What is the Lebesgue measure of $Q$?

(ii) Is the Hausdorff dimension of $Q$ equal to $\dim A$? The answer is yes in the quadratic case [Sh]

Mary Rees proved that the Lebesgue measure of $Q$ is positive [R] in the case when $A$ is the whole space of rational maps of degree $d$. On the other hand, Shishikura claims that in the quadratic family $z \mapsto z^2 + c$ the measure of the set of only finitely renormalizable points in $Q$ is equal to zero (here $Q$ is just the boundary of the Mandelbrot set). How do these results fit?

# On Invariant Measures for Iterations of Holomorphic Maps

## Feliks Przytycki

Let $U$ be an open subset of the Riemann sphere $\hat{\mathcal{C}}$. Consider any holomorphic mapping $f : U \to \hat{\mathcal{C}}$ such that $f(U) \supset U$ and $f : U \to f(U)$ is a proper map, (for a more general situation see [PS]). Consider any $z \in f(U)$. Let $z^1, z^2, ..., z^d$ be some of the $f$-preimages of $z$ in $U$ where $d \geq 2$. Consider curves $\gamma^i : [0,1] \to \hat{\mathcal{C}}$, $i = 1, ..., d$, also in $f(U)$, joining $z$ with $z^i$ respectively (i.e. $\gamma^i(0) = z, \gamma^i(1) = z^i$).

Let $\Sigma^d := \{1, ..., d\}^{\mathbf{Z}^+}$ denote the one-sided shift space and $\sigma$ the shift to the left, i.e. $\sigma((\alpha_n)) = (\alpha_{n+1})$. For every sequence $\alpha = (\alpha_n)_{n=0}^{\infty} \in \Sigma^d$ we define $\gamma_0(\alpha) := \gamma^{\alpha_0}$. Suppose that for some $n \geq 0$, for every $0 \leq m \leq n$, and all $\alpha \in \Sigma^d$, the curves $\gamma_m(\alpha)$ are already defined. Suppose that for $1 \leq m \leq n$ we have $f \circ \gamma_m(\alpha) = \gamma_{m-1}(\sigma(\alpha))$, and $\gamma_m(\alpha)(0) = \gamma_{m-1}(\alpha)(1)$.

Define the curves $\gamma_{n+1}(\alpha)$ so that the previous equalities hold (by taking $f$-preimages of curves already existing; if there are no critical values for iterations of $f$ in $\bigcup_{i=1}^{d} \gamma^i$ one has a unique choice). For every $\alpha \in \Sigma^n$ and $n \geq 0$ denote $z_n(\alpha) := \gamma_n(\alpha)(1)$.

The graph with the vertices $z$ and $z_n(\alpha)$ and edges $\gamma_n(\alpha)$ is called a *geometric coding tree* with the root at $z$. For every $\alpha \in \Sigma^d$ the subgraph composed of $z, z_n(\alpha)$ and $\gamma_n(\alpha)$ for all $n \geq 0$ is called a *geometric branch* and denoted by $b(\alpha)$. The branch $b(\alpha)$ is called *convergent* if the sequence $z_n(\alpha)$ is convergent in cl$U$. We define the *coding map* $z_{\infty} : \mathcal{D}(z_{\infty}) \to \text{cl}U$ by $z_{\infty}(\alpha) := \lim_{n \to \infty} z_n(\alpha)$ on the domain $\mathcal{D}(z_{\infty})$ of all such $\alpha$'s for which $b(\alpha)$ is convergent.

There are two basic examples:

1. $f : U \to U$ where $U$ is a simply-connected domain in $\hat{\mathcal{C}}$, $\deg f \geq 2$, and the iterates $f^n$ converge to a constant in $U$, in particular $U$ is an immediate basin of attraction of a sink for $f$ a rational map on $\hat{\mathcal{C}}$.

2. $U = \hat{\mathcal{C}}$, $f$ is a rational mapping.

It is known that except for a "thin" set in $\Sigma^d$ all branches are convergent (i.e. $\Sigma^d \setminus \mathcal{D}(z_{\infty})$ is "thin" and for every $x \in \text{cl}U$, the set $z_{\infty}^{-1}(x)$ is "thin"). These hold under very mild assumptions about the tree even allowing the existence of critical values in it. Proofs and a discussion of various possibilities of "thiness" can be found in [PS]. In particular one obtains the classical Beurling's Theorem that a holomorphic univalent function $R$ on the unit disc $I\!D$ has radial limits everywhere except on a set of logarithmic capacity zero, and for every limit point, the set in $\partial I\!D$ to which radii converge is also of logarithmic capacity 0. One just transports the map $z \mapsto z^2$ to $U := R(I\!D)$, and gets a type 1 situation. There is a 1-to-1 correspondence between the radii and geometric branches.

**General Problem.** How large is the image: $z_{\infty}(\mathcal{D}(z_{\infty}))$ ?

We shall specify this Problem separately in the basin of attraction case (the situation 1 above) and in the general situation.

To simplify the notation we have restricted ourselves to trees and codings from the full shift space. In the general situation it might be useful to consider also a topological

Markov chain, see [PS].

## THE CASE OF THE BASIN OF ATTRACTION

**Problem 1.1** If $f$ extends holomorphically to a neighbourhood of cl$U$, is every periodic point in $\partial U$ accessible from $U$ ?

**Comment.** Accessible means being $\varphi(1)$ for a continuous curve $\varphi : [0,1] \to \text{cl}U$ where $\varphi([0,1)) \subset U$ what is equivalent to being in the radial limit (i.e. $\lim_{r \nearrow 1} R(r\zeta)$ for $\zeta \in \partial I\!\!D$, $R$ denoting a univalent map from $I\!\!D$ onto $U$). For $g$ denoting the holomorphic extention of $R^{-1} \circ f \circ R$ to a neighbourhood of cl$I\!\!D$ and $\bar{R}$ the radial limit of $R$ wherever it exists, it is known that at every $g$-periodic $\zeta \in \partial I\!\!D$, $\bar{R}$ exists and $f$ at $\bar{R}(\zeta)$ is $f$-periodic (equivalently we could speak about $\sigma$-periodic points in $\Sigma^d$ and the mapping $z_\infty$, for a tree in $U$). Are there other periodic points in $\partial U$ ? It seems it does not matter if one assumes here that $f$ is defined only on a neighbourhood of $\partial U$. This is the case of an RB-domain $U$ (the boundary is repelling on the $U$ side) considered in [PUZ]. Problem 1.1 has a positive answer in the case where $f$ is a polynomial on $\mathbb{C}$ and $U$ is the basin of attraction to $\infty$, (Douady, Yoccoz, Eremenko, Levin), even if $U$ is not simply-connected, see [EL]. Here the fact $f^{-1}(U) \subset U$ helps.

**Problem 1.2.** In the situation of Problem 1.1 is every point $x \in \partial U$ of positive Lyapunov exponent (i.e. such that $\liminf_{n \to \infty} \frac{1}{n} \log |(f^n)'(x)| > 0$) accessible from $U$ ?

**Problem 1.3.** In the situation of Problem 1.1 is it true that the topological entropy $h_{\text{top}}(f|_{\partial U}) = \log \deg(f|_U)$ ?

**Comment** The $\geq$ inequality is known and easy. The problem is with the opposite one. It would be true if every point $x \in \partial U$ had at most $\deg(f|_U)$ pre-images in $\partial U$.

A positive answer to problem 1.2 would give a positive answer to 1.3. The reason is that topological entropy is approximated by measure-theoretic entropies for $f$-invariant measures which having positive entropies would have positive Lyapunov exponents (Ruelle's inequality). Then they would be images under $\bar{R}$ of $g$-invariant measures on $\partial I\!\!D$ which all have entropies upper bounded by $\log d$ (as $g$ is a degree $d$ expanding map on $\partial I\!\!D$).

**Problem 1.4.** Can there be periodic points or points with positive Lyapunov exponents in the boundary of a Siegel disc $S$ ? Is it always true that $h_{\text{top}}(f|_{\partial S}) = 0$?

## THE GENERAL CASE

We suppose here only that $f$ extends holomorphically to a neighbourhood of the closure of the limit set $\Lambda$ of a tree, $\Lambda = z_\infty \mathcal{D}(z_\infty)$. Then $\Lambda$ is called a quasi-repeller, see [PUZ]. Denote the space of all probability $f$-invariant ergodic measures on the closure of a quasi-repeller $\Lambda$ by $M(\Lambda)$. The space of measures in $M(\Lambda)$ which have positive entropy will be denoted by $M^+(\Lambda)$.

**Problem 2.1.** Is it true that every $m \in M(\Lambda)$ is the image of a measure on the shift space $\Sigma^d$ through a geometric coding tree with $z$ in a neighbourhood of cl$\Lambda$. What about measures in $M^+(\Lambda)$ ? The same questions for $f$ a rational mapping of degree $d$ on $U = \hat{\mathbb{C}}$ and measures on the Julia set $J(f)$.

**Comment.** It is easy to see at least, due to the topological exactness of $f$ on the Julia set $J(f)$ (for every open $V$ in $J(f)$ there exists $n > 0$ so that $f^n(V) = J(f)$), that for every $z$ except at most two, $z_\infty(\mathcal{D}(z_\infty)$ is dense in $J(f)$. The answer is of course positive in the case $f$ is expanding on $\Lambda$ because then $z_\infty$ is well defined and continuous on $\Sigma^d$, hence $\Lambda$ is closed.

**Problem 2.2** For which $m \in M^+(\Lambda)$ for every "reasonable" function $\varphi : \Lambda \to \mathbb{R} \cup \pm\infty$ (for example Hölder, into $\mathbb{R}$ or allowing isolated values $-\infty$ with $\exp\varphi$ nonflat there, as $\log|g|$, $g$ holomorphic) do the probability laws like Almost Sure Invariance Principle, Law of Iterated Logarithm, or Central Limit Theorem hold for the sequence of sums $S_n(\varphi) = \sum_{j=0}^{n-1} t_j$ of the random variables $t_j := \varphi \circ f^j - \int \varphi \, dm$ provided $\sigma^2(\varphi) = \lim \frac{1}{n} \int S_n(f)^2 dm > 0$ ?

**Comment.** If the measure is a $z_\infty$-image of a measure on $\Sigma^d$ with a Hölder continuous Jacobian (a Gibbs measure for a Hölder continuous function) then the probability laws hold, see [PUZ]. The positive answer in Problem 2.1 would be very helpful in solving Problem 2.2.

The class of measures for which Problem 2.2 has not been solved, but does not seem out of reach, are equilibrium states for Hölder continuous functions, say on the Julia set in the case $f$ is rational. In this case the transfer (Ruelle-Perron-Frobenius) operator is already understood to some extent [DU], [P]. A proof seems to depend on finding an appropriate space of functions on which the maximal eigenvalue has modulus strictly larger than supremum over the rest of the spectrum (by the analogy to the expanding case, [Bowen]).

Actually these equilibrium states are $z_\infty$-images of measures on $\Sigma^d$. The Jacobians of these equilibrium states have modulus of continuity bounded by Const$(m)(\log(1/t))^{-m}$ for any $m > 0$ (I don't know if it is Hölder). The Jacobian of the pull-back of the equilibrium measure to $\Sigma^d$ is not wild. This gives a chance to prove that mixing in $\Sigma^d$ is polynomially fast.

**Problem 2.3** Is it true for every $m \in M^+(\Lambda)$ that $m$ is absolutely continuous with respect to $H_\kappa$ (where $H_\kappa$ is the Hausdorff measure in dimension $\kappa = \mathrm{HD}(m)$) iff $HD(m) = HD(\mathrm{cl}\Lambda)$?

**Comment.** In such a generality I would expect a negative answer. One should probably restrict the family of measures under consideration and/or impose additional assumptions on the mapping $f$.

If $f$ is expanding on $\Lambda$ then the answer is positive for all measures in $M^+(\Lambda)$ with Hölder continuous Jacobian. This is basically Bowen's theorem.

In the discussion here we assume that on every set $E$ on which $f$ is 1-to-1 the measure $(f|_E)^{-1}(m)$ is equivalent to $m$, and we write $\mathrm{Jac}_m f(z) = \frac{d(f|_E)^{-1}(m)}{dm}(z)$.

When the Jacobian exists in this sense we can replace the absolute continuity hypothesis $m \ll H_\kappa$ or the alternative singularity hypotesis $m \perp H_\kappa$ with another pair of alternative hypotheses.

**Problem 2.4.** In what class of measures in $M^+(\Lambda)$ does the property: the family $S_n(\log \mathrm{Jac}_m(f) - \kappa \log|f'|)$ is not uniformly bounded in $L^2(m)$, imply $m \perp H_\kappa$ and $HD(m) < HD(\mathrm{cl}\Lambda)$.

**Comment.** The answer is positive for $f$ expanding and Jacobian Hölder continuous.

It is positive also if $m = z_\infty(\mu)$ for any Gibbs measure $\mu$ for a Hölder continuous function on $\Sigma^d$. The singularity $\perp$ follows then from the positive answer to Problem 2.2 in this special case, see [PUZ]. From the probability laws one can deduce a stronger singularity, for example with respect to the measure $H_\Phi(\kappa, c)$ which is the Hausdorff measure for the function

$$\Phi(\kappa, c)(t) = t^\kappa \exp c \sqrt{\log \frac{1}{t} \log \log \log \frac{1}{t}}$$

for all

$$c < \sqrt{2\sigma^2(\log \mathrm{Jac}_\mu(s) - \kappa \log|f'| \circ z_\infty) / \int \log|f'| dm}.$$

The inequality $HD(m) < HD(\mathrm{cl}\Lambda)$ follows from [Z1].

**Problem 2.5.** In what class of measures in $M^+(\Lambda)$ s does the property: the family $S_n(\log \mathrm{Jac}_m(f) - \kappa \log|f'|)$ is uniformly bounded in $L^2(m)$, imply $m \ll H_\kappa$ ?

**Comment.** Again the answer is positive for $f$ expanding and Jacobian Hölder continuous.

If $m = z_\infty(\mu)$ then the boundness of the family $S_n(\varphi)$ where $\varphi := \log \mathrm{Jac}_\mu(f) - \kappa \log|f'| \circ z_\infty$ occurs precisely when $\sigma^2(\log \mathrm{Jac}_\mu(f) - \kappa \log|f'| \circ z_\infty) = 0$ assuming the series $\sum_{n=1}^\infty n \int |\varphi \cdot (\varphi \circ s^n)| d\mu$ is convergent. This is equivalent to the existence of a function $u$ in $L^2(\mu)$ so that $\varphi = u \circ s - u$. Then we say that we can solve the cohomology equation for $\varphi$. Then we can also solve the cohomology equation for $\log \mathrm{Jac}_m(f) - \kappa \log|f'|$ on $\Lambda$. The naive way to compare $m$ with $H_\kappa$ is to prove that the sequence $S_n(\log \mathrm{Jac}_m(f) - \kappa \log|f'|)(z)$ is bounded at almost every $z \in \Lambda$. In the expanding case this allows comparison of the $m$-measure and the radius to the $\kappa$ power of little discs, so the naive method happens to be successful. In the general case we do not have even pointwise boundness, because the function $u$ is only in $L^2(\mu)$.

The problem has the positive answer in the following special cases:

1. In the RB-domain case, where m is equivalent to a harmonic measure on the boundary of a simply-connected domain $U$, see [PUZ] and [Z2]. Then $m = \bar{R}(\mu)$ where $\mu$ is equivalent to the Lebesgue measure on $\partial I\!\!D$. $\log|R'|$ happens to be within a bounded distance from any harmonic extension of $u$ to a neighbourhood of $\partial I\!\!D$, in particular radial limits for $\log|R'|$ exist a.e.. In [Z2] it is proved in fact that all this implies that $\partial U$ is analytic, giving the answer to Problem 2.3 in this case.

2. In the case where $f$ is a rational map on $\hat{\mathbb{C}}$ and $m$ is a measure with maximal entropy (in which case Jacobian$\equiv \deg f$). Then again a careful look at $u$ proves that $f$ is either $z \mapsto z^n$ or is a Tchebysheff polynomial (in respective holomorphic coordinates on $\hat{\mathbb{C}}$) or else $J(f) = \hat{\mathbb{C}}$ and $f$ has a parabolic orbifold, see [Z1].

In the general case it seems hopeful to treat any harmonic extension of $u$ as a logarithm of a derivative of a "Riemann mapping". In the case $m = z_\infty(\mu)$ one can average $u$ over cylinders in $\Sigma^d$ extending $u$ to the vertices $z_n(\alpha)$ of the tree.

The mapping $z_\infty$ can be viewed as a dynamical version of a Riemann maping. We can formulate the following problem:

**Problem 2.6.** Which theorems about the boundary behaviour of Riemann maps hold for geometric coding trees?

**Comment.** Beurling Theorems hold, see the discussion in Section 1.

One has a natural dictionary:

| | |
|---|---|
| For $R$: | For $z_\infty$: |
| prime end | a geometric branch |
| impression | $I(\alpha) = \cap_{n=0}^{\infty} z_\infty\{\beta : \beta_i = \alpha_i, i = 0, ..., n\}$ |
| the set of principal points | the limit set for the vertices $z_n(\alpha)$ of $b(\alpha)$. |

**Problem 2.7.** Is it true that $\sup_{m \in M^+(\Lambda)} \mathrm{HD}(m) = \mathrm{HD}(\mathrm{cl}\Lambda)$? Does $\sup_{m \in M(\Lambda)}$ help?

**Comment.** Of course a negative answer to this Problem for some $\Lambda$ and positive to Problem 2.3 would mean that $m \perp H_{\mathrm{HD}(m)}$ for all $m$.

Problem 2.7 has positive answer in the expanding and subexpanding cases where sup is attained, it is so even for a positive measure set of rational mappings on $\hat{\mathbb{C}}$ for which absolutely continuous invariant measures exist (with respect to the Lebesgue), see [R]. The problem has also a positive answer for rational mappings with neutral points but without critical points in the Julia set. But then it may happen that supremum is not attained, see [ADU] and [L].

# Section 4: Iterates of Entire Functions

## Open Questions in Non-Rational Complex Dynamics
### Robert Devaney

The dynamics of complex analytic functions have been studied by many authors during the past decade. Much of this work has been confined to the study of either rational or polynomial maps. The study of other analytic functions is still in its infancy and there are many unsolved problems in this area. In this note we describe a few of these problems.

**1. Entire functions.** The dynamics of entire functions are quite different from the dynamics of rational maps, mainly because of the essential singularity at infinity. By the Picard theorem, any neighborhood of this singularity is mapped infinitely often over the entire plane missing at most one point. This injects considerable hyperbolicity into the map and often causes the topology of the Julia set of the map to be vastly different from that of a rational map. In addition, the No Wandering Domains Theorem of Sullivan does not hold for this class of maps, so there may be both wandering domains and domains at infinity in the stable sets.

There is one class of entire maps whose dynamics are fairly well understood, namely the entire maps that have finitely many asymptotic and critical values (maps of finite type). With few exceptions (notably examples of Baker [B], Herman [H], and Eremenko and Lyubich [EL], most work has centered around this class of maps. Extending the study to a wider class of maps is an important problem.

*Problem:* Find a collection of representative examples of entire maps whose dynamics may be understood.

As a starting point, one might ask

*Problem:* What are the dynamics of maps of the form $\lambda e^z \sin z$ or $\lambda e^z \cos z$?

**2. Entire functions of finite type.** Most of the work thus far on the dynamics of entire maps has been concentrated on the class of finite type maps. These are the maps which have only finitely many singular (i.e., critical and asymptotic) values. This class includes $\lambda e^z$, $\lambda \sin z$, and $\lambda \cos z$. It is known [GK,EL] that the No Wandering Domains theorem holds for this class, and that the Julia sets of these maps often contain Cantor bouquets [DT].

**3. The exponential map.** Of all entire maps, the exponential family $E_\lambda(z) = \lambda e^z$ has received the most attention. This is natural since $E_\lambda$, like the well-studied quadratic family $Q_c(z) = z^2 + c$, has only one singular value, the asymptotic value at 0. Thus this family is a "natural" one parameter family.

The parameter space for $E_\lambda$ has been studied in [DGK]. However, there remain significant gaps in this picture. It is known that there exists Cantor sets of curves (called hairs) in the parameter plane for which the corresponding exponential maps have Julia sets that are the whole plane.

*Problem:* Describe completely the set of $\lambda$-values for which the Julia set of $E_\lambda$ is **C**.

*Problem:* Many of these $\lambda$-values lie on curves or hairs. Are these hairs $C^\infty$? Analytic? Where and how do they terminate?

There are some interesting topological structures embedded in the dynamics of the exponential that warrant further study. For example, it is known that for $\lambda > 1/e$, $J(E_\lambda) =$ **C**. However, if $\lambda, \mu > 1/e$, then $E_\lambda$ and $E_\mu$ are not topologically conjugate [DG]. If one looks at the invariant set consisting of $\{z \| 0 \le \operatorname{Im} E_\lambda^n(z) \le \pi \text{ for all } n\}$, it is known that this set is a Knaster-like continuum.

*Problem:* Are each of these Knaster-like continua homeomorphic? (for any $\lambda, \mu > 1/e$)

**4. The Trigonometric Functions**. The parameter spaces for families such as $S_\lambda(z) = \lambda \sin z$ or $C_\lambda(z) = \lambda \cos z$ also deserve special attention. They also contain curves on which the Julia set is the entire plane. The fundamental difference here is that $C_\lambda$ and $S_\lambda$ have no finite asymptotic values (only critical values), whereas the opposite is true for $E_\lambda$.

*Problem:* Describe the structure of the parameter space for $C_\lambda$ and $S_\lambda$.

One fundamental difference between the trigonometric and exponential families is the following. Both maps are known to possess Cantor bouquets [DT] in their Julia sets. And any two planar Cantor bouquets are homeomorphic [AO]. Finally, McMullen [Mc] has shown that these Cantor bouquets always have Hausdorff dimension 2. However, the Lebesgue measure of these bouquets is quite different: they always have measure zero in the exponential case, but infinite measure in the trigonometric case.

*Problem:* What is the measure and dimension of the hairs i the parameter space for $E_\lambda, S_\lambda$ and $C_\lambda$.

**5. Other families of non-rational maps.** Newton's method applied to non-rational maps offers a fertile area for further investigation. Outside of the work of Haruta [Ha] and van Haesler and Kriete [HK ], there is little that is known. So a general problem is:

*Problem:* Describe the dynamics of Newton's method applied to general classes of entire functions?

This, of course, immediately leads to the question of iteration of meromorphic functions. Some work has been done here in case the map has polynomial Schwarzian derivative [DK] or when the map has finitely many singular values [BK]. But not much else is known.

Finally, there is an intriguing object called the tricorn introduced by Milnor [M] as one of his basic slices of parameter space for higher dimensional maps. This object arises as the analogue of the Mandelbrot set for the anti-holomorphic family $A_c(z) = \overline{z}^2 + c$. It is known [La] that the tricorn is not locally connected, but it also contains smooth arcs in the boundary (with no decorations attached) [W]. As this object arises in slices of the cubic connected locus, it certainly warrants further study. Winters [W] also has introduced a family of fourth degree polynomials whose parameter space is "naturally" $\mathbf{R}^3$ and which contains perpendicular slices given by the Mandelbrot set and the tricorn. Winters suggests that this family can model cubics since there are only two critical orbits.

# Wandering Domains for Holomorphic Maps

## A. Eremenko and M. Lyubich

Let $f$ be a rational or entire function. A connected component $D$ of the complement of the Julia set $J(f)$ is called wandering domain if for all $m > n \geq 0$ we have $f^m D \cap f^n D = \emptyset$, where $f^m$ stands for $m$-th iterate of $f$. One of the most important theorems in holomorphic dynamics due to D. Sullivan states that rational functions have no wandering domains [11]. We ask for possible generalizations of this theorem. All known proofs of Sullivan's theorem use heavily the fact that the space of quasiconformal deformations of a rational function is finitely dimensional (see e.g. [3]).

Here is one situation where a similar result could be proved. We say that an entire function $f$ belongs to the class $S$ if there is a finite set of points $\{a_1, \ldots, a_q\}$ such that

$$f : \mathbf{C} \backslash f^{-1}\{a_1, \ldots, a_q\} \to \mathbf{C} \backslash \{a_1, \ldots, a_q\}$$

is a covering map. The space of quasiconformal deformations of an entire function of the class $S$ has finite dimension and the following result can be proved by extending the Sullivan's method: entire functions of the class $S$ have no wandering domains [6], [8], [2].

On the other hand it is known that wandering domains $D$ may exist for some entire functions $f$. The examples with the following properties have been constructed:

1). $f^n D \to \infty$, [1], [2], [5], [9]. The example in [5] has an additional property that the iterates $f^n$ are univalent in $D$.

2). The orbit $\{f^n D\}$ has infinitely many limit points, including $\infty$, [5].

**Question 1** *Does there exist an entire function $f$ with a wandering domain $D$ such that the orbit $\{f^n D\}$ is bounded?*

Remark that there are entire functions not in the class $S$, for which the negative answer can be obtained easily. We say that the function $f$ has order less then one half if

$$\log \log^+ |f(z)| \leq \alpha \log |z|, \quad |z| > r_0$$

for some $\alpha < 1/2$. It follows from a classical theorem by Wiman and Valiron (see, for example, [10]) that such functions have the following property: there exists a sequence $r_k \to \infty$ such that

$$|f(r_k e^{i\theta})| > r_k, \quad 0 \leq \theta \leq 2\pi.$$

It follows that there is an increasing sequence of domains $G_k$, $\cup G_k = \mathbf{C}$ such that the restrictions of $f$ on $G_k$ are polynomial-like maps [4]. So $f$ has no wandering domains with bounded orbit because polynomial-like maps have no wandering domains.

Now we consider a special type of wandering domains whose orbits tend to a finite point $z_0$. Let $\varphi$ be a germ of holomorphic function with the point $z_0$ fixed. Suppose that $\lambda = \varphi'(z_0) = \exp 2\pi i \alpha$, $\alpha$ irrational. It was proved by Fatou [7] that in this situation $\varphi^n(z)$ cannot tend to $z_0$ in an *invariant* domain. So we have the following

**Question 2** *Is it possible that $\varphi^n(z) \to z_0$ uniformly in some domain $D$?*

In the case when $\varphi$ can be analytically continued to an entire function positive answer would imply the existence of wandering domain whose orbit tends to $z_0$. It would be also interesting to know the answer to the question 2 with other additional assumptions on the germ $\varphi$, for example, when $\varphi$ is a germ of an algebraic function.

Finally remark that the answer to the following question is also unknown

**Question 3** *Under the assumptions of Question 2 can it happen that there is an orbit tending to $z_0$?*

*

# Section 5: Newton's Method

## Bad Polynomials for Newton's Method
### Scott Sutherland

Newton's method for solving $f(z) = 0$ corresponds to iteration of $z \mapsto z - f(z)/f'$ which is a degree $d$ rational map of $\overline{\mathbb{C}}$ in the case where $f$ is a polynomial of degree with distinct roots. Newton's method has long been an important source of examples theorems in complex dynamical systems (for example, the work of Schröder [Sch, Sc] Fatou [Fa], and more recently Douady and Hubbard [DH]), as well as being one of most commonly used numerical schemes for approximating roots. See [HP] and [Sm] for introduction to the dynamics of Newton's Method.

Describing the set of polynomials for which the corresponding Newton's method periodic sinks which are not roots is an important open problem, (problem 6 of [Sm]). shall refer to such polynomials as "bad polynomials". This question is essentially answe for cubic polynomials by the work of Tan Lei [Ta] and Janet Head [He], in which the m comprehensive task of giving a combinatorial description of the parameter space for Newto method is undertaken. A complete description of the parameter space for higher degrees seems some way off, however.

In order to answer Smale's question for higher degree polynomials, it may be helpful consider the relationship between the "relaxed Newton's method"

$$N_{h,f}(z) = z - h \frac{f(z)}{f'(z)}$$

and the "Newton Flow" $\mathcal{N}_f$ given by the ordinary differential equation

$$\dot{z} = -\frac{f(z)}{f'(z)}.$$

One sees immediately that the map is an Euler approximation to the flow using step $h$. The attractors of $\mathcal{N}_f$ are sinks located at the zeros of $f(z)$, $\infty$ is the only source, the other fixed points are at the singularities corresponding to the critical points of $f$. can rescale time for $\mathcal{N}_f$ to obtain $\dot{z} = -f(z)\overline{f'(z)}$ (or alternatively $\dot{z} = -\nabla\|f(z)\|^2$), fr which we can easily see that these singularities are hyperbolic saddles. Furthermore, solut curves of $\mathcal{N}_f$ are mapped by $f$ to straight lines emanating from the origin. Thus, if $f$ two critical values with the same argument, then the flow $\mathcal{N}_f$ is degenerate in the sense t there are solution curves which begin at one singularity and terminate at another. Refe [JJT], [Sa], [STW], [Sm], and [Su] for more details about $\mathcal{N}_f$.

**Conjecture 1.** Let $f_1$ be a bad polynomial of degree $d$, that is one for which Newton method has an attractor which is not a root of $f$. Then there is a one-parameter family polynomials $\{f_h\}_{0 < h \leq 1}$ which are bad for the relaxed Newton's method $N_{h,f_h}$. Furthermore as $h \to 0$, the corresponding flow $\mathcal{N}_{f_h}$ tends to a flow $\mathcal{N}_{f_0}$ which is degenerate.

This conjecture is consistent with the following, as explained below.

**Conjecture 2.** Let $f$ be a polynomial of degree $d$ with all its roots in the unit disk, $\alpha$ be a root of multiplicity $m$ for $f$, and let $A_h^*(\alpha)$ be the immediate attractive basin of $\alpha$ the map $N_{h,f}$. Then the intersection of the set

$$\mathcal{A} = \bigcap_{0 < h \leq m} A_h^*(\alpha)$$

with any circle of radius $R \geq 3$ contains arcs whose total length is at least $\frac{2\pi R}{cd}$, where $c$ constant not depending on $\alpha$, $f$, or $d$.

This second conjecture says that there is a definite neighborhood of the singular jectories of $\mathcal{N}_f$ in which the Julia set of $N_{h,f}$ must be contained for all $h \in (0, m]$. Si the periodic orbits for $N_{h,f}$ which are not roots must be contained in the complement $\bigcup_{f(\alpha)=0} A_h^*(\alpha)$, conjectures 1 and 2 taken together give some idea of the structure of parameter space for $N_{h,f}$.

Conjecture 2 has been partially established by Benzinger [Be] (for all $h$ sufficiently n 0), and is a generalization of the main result of [Su], which shows this for $h = 1$. I beli that with slight modifications, the proof in [Su] can be made to work for $0 < h \leq m$, wh should nearly complete the proof of conjecture 2.